\documentclass[12pt]{article}
\usepackage{amsmath}
\usepackage{amscd}
\usepackage{mathrsfs}
\usepackage{amsfonts}
\usepackage{amsmath,amssymb}
\baselineskip 5.5pt \lineskip 5.5pt \numberwithin{equation}{section}

\def \Z{\hbox{$Z\hskip -5.2pt Z$}}

\def \C{\hbox{$C\hskip -5pt \vrule height 6pt depth 0pt \hskip 6pt$}}

\def\qed{\ \ \ifhmode\unskip\nobreak\fi\ifmmode\ifinner
         \else\hskip5pt\fi\fi
 \hbox{\hskip5pt\vrule width4pt height6pt depth1.5pt\hskip 1 pt}}
\def\a{\alpha}

\def\b{\beta}
\def\d{\delta}

\def\g{\gamma}

\def\L{\mathcal{L}}

\def\cl{\centerline}

\def\vs{\vspace*}

\def\C{\mathbb{C}}

\def\Z{\mathbb{Z}}

\newtheorem{theo}{Theorem}[section]
\newtheorem{lemm}[theo]{Lemma}

\newtheorem{defi}[theo]{Definition}
\newtheorem{prop}[theo]{Proposition}

\begin{document}
\cl {\large\bf
 \vs{10pt} Irreducible representations of the twisted Heisenberg-Virasoro}
\cl {{\large\bf algebra of Hom type}
\noindent\footnote{Supported by the National Science Foundation of
China (Nos. 11047030 and 11771122).
}} \vs{6pt}

\cl{Qiuli Fan, Yongsheng Cheng}

\cl{ \small School
of Mathematics and Statistics, Henan
University, Kaifeng 475004, China} \vs{6pt}

\vs{6pt}
{\small
\parskip .005 truein
\baselineskip 10pt \lineskip 10pt
\noindent{{\bf Abstract:}\, Hom-Lie algebras are non-associative algebras generalizing Lie algebras by twisting the
Jacobi identity by a linear map. In this paper, we mainly study the irreducible representation of the twisted Heisenberg-Virasoro algebra of Hom-type, which can be induced by the irreducible representations of its induced Lie algebra. In particular, we construct some kinds of irreducible representations of the twisted Heisenberg-Virasoro algebra of Hom-type. }
\vs{5pt}

}

\noindent{\bf MR(2020) Subject Classification} 17A30, 17B61, 17B68.
\parskip .001 truein\baselineskip 8pt \lineskip 8pt

\vs{6pt}
\par
\cl{\bf\S1. \ Introduction}
\setcounter{section}{1}\setcounter{theo}{0}\setcounter{equation}{0}

The twisted Heisenberg-Virasoro Lie algebra at level zero was first introduced in \cite{[ACKP]}, which is the
universal central extension of the Lie algebra of differential
operators on a circle of order at most one.
 This Lie algebra has an infinite dimension Heisenberg subalgebra and
a Virasoro subalgebra. These subalgebras, however, do not form a
semidirect product, but instead, the natural action of the Virasoro
subalgebra on the Heisenberg subalgebra is the twisted with a
2-cocycle.

The motivations to study Hom-Lie structures are related to physics
and to deformations of Lie algebras, in particular Lie algebras of
vector fields. As a generalization of Lie algebras, Hom-Lie algebras were
introduced by Hartwig, Larsson and Silvestrov in \cite{[HLS]} as part of a
study of deformations of the Witt and the Virasoro algebras.
A Hom-Lie algebra is a triple $(\L,[\cdot,\cdot]_{\varphi},\varphi)$ in which the bracket satisfies skew-symmetry and a twisted Jacobi identity along the linear map $\varphi$ which is called the Hom-Jacobi identity. The main feature of these algebras
is that the identities defining the structures are twisted by homomorphisms.
The paradigmatic examples are $q$-deformations of
Witt and Virasoro algebras, Heisenberg-Virasoro algebra and other algebraic structure
constructed in pioneering works
\cite{ [CS2], [DSS], [JL], [XJL]}.

The representation theory of Hom-Lie algebras
was introduced by Sheng in \cite{[S]}, in which, Hom-cochain complexes, derivation,
central extension, derivation extension, trivial representation and adjoint representation
of Hom-Lie algebras were studied. In addition,
(co)homology and deformations theory of (Bi)-Hom-algebras were studied in
\cite{[CQ], [CS1], [CY], [LCC]}. The diversity of the twist map of
g makes this topic interesting and complicated.

In nowadays, one of the most pervasive trends in mathematics has to do with representation theory. This topic is an important tool in most parts of Mathematics and Physics. Yunhe Sheng introduced representations of multiplicative hom-Lie
algebras in \cite{[S]}. The relationship between representations of  Hom-Lie algebra and its induced Lie algebra is introduced in \cite{[ABM], [L]}. It can be obtained that the irreducible representation of Lie algebra of Hom type by the irreducible representation of semi-simple Lie algebra when $\phi\in End(V)$ is invertible from \cite{[ABM]}, and according to this theory, Hom-type $sl(2)$-modules are divided into four types which include \\(1) finite-dimensional irreducible modules,\\(2) irreducible infinite-dimensional lowest weight Hom-sl(2)-modules, \\(3) irreducible infinite-dimensional highest weight Hom-sl(2)-modules, \\(4) irreducible infinite-dimensional Hom-sl(2)-module (the Hom-module of intermediate series).\\ Based on this, we aim to classify irreducible representation over the twisted Heisenberg-Virasoro algebra of Hom type in this paper.

Hom-Lie algebras are non-associative algebras generalizing Lie algebras by twisting the
Jacobi identity by a linear map. In \cite{[LTL]}, the authors determined some kinds of nontrivial multiplicative Hom-Lie structures on the twisted Heisenberg-Virasoro algebra. In this paper, we mainly study the irreducible representation of the twisted Heisenberg-Virasoro algebra of Hom-type, which can be given by the irreducible representations of its induced Lie algebra. In particular, we construct some kinds of irreducible representations of the twisted Heisenberg-Virasoro algebra of Hom-type.
The purpose of this paper is to study representations of Bihom-Lie algebras.
The paper is organized as follows. In Section 2 we give the basic definitions of twisted Heisenberg-Virasoro algebras and their some Harish-Chandra modules. Then we review definitions of representation and introduce the relationship between representations of  Hom-Lie algebras and representations of their induced Lie algebras. In section 3, we get the irreducible representations over twisted Heisenberg-Virasoro of Hom type with respect to Harish-Chandra modules of the series. We need to note these irreducible representations are not always weight modules.
\vs{6pt}

\cl{\bf\S2.\  Notations and Preliminaries}
\setcounter{section}{2}\setcounter{theo}{0}\setcounter{equation}{0}
In this section we provide some notations, definitions and preliminary results.
\begin{defi}
\label{heivir}
The twisted Heisenberg-Virasoro algebra $\L$ is defined to be the Lie algebra with a basis $$\{L_n,I_n,C_L,C_LI,C_I|n\in\Z\}$$
subject to the following commutation relations
$$\begin{array}{cc}
[L_n,L_m]=(n-m)L_{n+m}+\d_{n,-m}\frac{n^3-n}{12}C_L,\\[6pt]
[L_n,I_m]=-mI_{n+m}+\d_{n,-m}(n^2+n)C_{LI},\\[6pt]
[I_n,I_m]=n\d_{n,-m}C_I,\\[6pt]
[\L,C_L]=[\L,C_{LI}]=[\L,C_I]=0.
\end{array}$$
\end{defi}

The following theorem coming from \cite{[LTL]}.
\begin{theo}
\label{nonendHV}
Every nonzero endomorphism of the twisted Heisenberg-Virasoro algebra $\varphi$ has the form as follows:
\begin{align}
&\varphi(L_n)=\frac{1}{k}a^nL_{kn}+a^n(cn+d)I_{kn}+\delta_{n,0}\frac{1-k^2}{24k}C_L
+\delta_{n,0}(ck+d)C_{LI}+\delta_{n,0}\frac{k}{2}(c^2-d^2)C_I,\\
&\varphi(I_n)=a^nbI_{kn}+\delta_{n,0}b((1-k)C_{LI}-(kc+kd)C_I),\\
&\varphi(C_L)=kC_L-24kcC_{LI}-12kc^2C_I,\\
&\varphi(C_{LI})=kbC_{LI}+kbcC_I,\\
&\varphi(C_I)=kb^2C_{I},
\end{align}
where $k\in\Z^*$, $a,b,c,d\in\C$, $n\in\Z$ and $a\neq0$.
\end{theo}

Liu and Jiang in 2008 gave a classification of the twisted Heisenberg-Virasoro algebra modules of the intermediate
series in \cite{[LJ]}. There are seven families of indecomposable modules of the intermediate series over the
the twisted Heisenberg-Virasoro algebra, that is, nontrivial indecomposable weight modules with weight multiplicity is at most one. They are Virasoro modules without central charge.
\begin{theo}
\label{hcmod}
Any Harish-Chandra module of intermediate series over the Heisenberg-Virasoro
algebra is isomorphic to one of the following modules (if they are irreducible) or their nontrivial simple
subquotients.
\begin{align}\label{a1}
&A(\a,\b,F): L_n(v_t)=(\a+t+{\b}n)v_{n+t},
I_n(v_t)=Fv_{n+t};\\
\nonumber &A(\a,F): L_n(v_t)=(t+n)v_{n+t}, t\neq0, L_n(v_0)=n(n+\a)v_n, \\
\label{a2} &\ \ \ \ \ \ \ I_n(v_t)=0, t\neq0, I_n(v_0)=nFv_n;\\
\nonumber &B(\a,F): L_n(v_t)=tv_{n+t}, t\neq{-n}, L_n(v_{-n})=-n(n+\a)v_0, \\
\label{a3} & \ \ \ \ \ \ \ I_n(v_t)=0, t\neq{-n}, I_n(v_{-n})=nFv_0;\\
\label{a4} &U_F: L_n(v_t)=tv_{n+t}, I_n(v_t)=0, t\neq{-n}, I_n(v_{-n})=nFv_0;\\
\label{a5} &V_F: L_n(v_t)=(t+n)v_{n+t}, I_n(v_t)=0, t\neq0, I_n(v_0)=nFv_n;\\
\label{a6} &\widetilde{U}_F: L_n(v_t)=tv_{n+t}, t\neq{-n},
L_n(v_{-n})=0, I_n(v_t)=0, t\neq{-n}, I_n(v_{-n})=nFv_0;\\
\nonumber &\widetilde{V}_F: L_n(v_t)=tv_{n+t}, t\neq{-n},
L_n(v_{-n})=0, \\
\label{a7}  &\ \ \ \ \ I_n(v_t)=0, t\neq0\ or\ n=t=0,  I_n(v_0)=Fv_n,
\end{align}
where $\a,\b,\g,F\in\C$, $n,t\in\Z$ and $C_L, C_I, C_{LI}$ act trivial.
\end{theo}

It is easy to know that the above modules are not isomorphic each other. Let
$v_{n}, n\in \mathbb{Z}$ be a basis of the Harish-Chandra modules of the intermediate series.

\begin{defi}\label{defhomlie}
A Hom Lie algebra is triple $(\L,[\cdot,\cdot]_{\varphi},\varphi)$ consisting of a vector space $\L$, a bilinear map $[\cdot,\cdot]_{\varphi}:\L\times\L\rightarrow\L$ and a linear map $\varphi:\L\rightarrow\L$ that satisfies, for any $x, y, z\in \L$,
\begin{center}
$[x,y]_{\varphi}=-[y,x]_{\varphi}$, \\
$[\varphi(x),[y,z]_{\varphi}]_{\varphi}+[\varphi(y),[z,x]_{\varphi}]_{\varphi}+[\varphi(z),[x,y]_{\varphi}]_{\varphi}=0$.
\end{center}
A Hom Lie algebra $(\L,[\cdot,\cdot]_{\varphi},\varphi)$ is called a multiplicative Hom-Lie algebra, if
$$\varphi([x, y]_{\varphi})=[\varphi(x), \varphi(y)]_{\varphi}). $$
A multiplicative Hom-Lie algebraIt is said regular if $\alpha$ is invertible.
\end{defi}

The next proposition coming from \cite{[ABM]} told us that, using Lie algebra and its homomorphism, we can get the corresponding
multiplicative Hom-Lie algebra.
\begin{prop}
\label{yautwist}
Let $(\L,[\cdot,\cdot])$ be a Lie algebra and $\varphi$ be a Lie algebra homomorphism. Then $(\L_{\varphi},[\cdot,\cdot]_{\varphi}:=\varphi\circ [\cdot,\cdot], \varphi)$ is a multiplicative Hom-Lie algebra.
\end{prop}

We call $(\L,[\cdot,\cdot])$ the induced Lie algebra of the multiplicative Hom-Lie algebra $(\L_{\varphi},[\cdot,\cdot]_{\varphi}, \varphi)$. By Proposition \ref{yautwist}, we can easily obtain next proposition.
\begin{prop}
\label{indliealg}
Let $(\L_{\varphi},[\cdot,\cdot]_{\varphi}, \varphi)$ be a regular Hom-Lie algebra, then $(\L, \varphi^{-1}\circ[\cdot,\cdot])$ is the induced Lie algebra of $(\L_{\varphi},[\cdot,\cdot]_{\varphi}, \varphi)$.
\end{prop}
\begin{defi}\label{rephom-lie}
A representation of a Hom Lie algebra $\L$ on the vector
space $V$ with respect to $\phi\in End(V)$ is a linear map $\rho_{\phi}:\L\rightarrow End(V)$, such that for any $x$, $y\in\L$, $v\in V$, satisfies:
\begin{equation}
\label{homrep1}
\rho_{\phi}([x,y]_{\varphi})\circ\phi(v)=\rho_{\phi}(\varphi(x))\rho_{\phi}(y)(v)-\rho_{\phi}(\varphi(y))\rho_{\phi}(x)(v).
\end{equation}
We call $(V,\rho_{\phi},\phi)$ a $\L$-module.

For a multiplicative Hom-Lie algebra $\L$, if $(V,\rho_{\phi},\phi)$ satisfies (\ref{homrep1}) and
\begin{equation}
\label{homrep2}
\rho_{\phi}(\varphi(x))\circ\phi(v)=\phi\circ\rho_{\phi}(x)(v),
 \end{equation}
then we call $(V,\rho_{\phi},\phi)$ a representation of a multiplicative Hom-Lie algebra $\L$.
\end{defi}
\begin{defi}\label{subrephom-lie}
Let $(V, \rho_{\phi},\phi)$ a representation of a Hom-Lie algebra $\L$. If $W$ is an invariant subspace of $V$ under $\rho_{\phi}(x)$ for any $x\in \L$, we call $(W, \rho_{\phi},\phi|_{W})$ a subrepresentation of $\L$.

The representation $(V,\rho_{\phi},\phi)$ of $\L$ is called irreducible, if it has nontrivial $\L$-subrepresentation.
\end{defi}

In \cite{[L]}, the author introduced the relationship between representations of the multiplicative Hom-Lie algebras and representations of their induced Lie algebras.
\begin{theo}\label{rep-homlielie}
Let $(\L,[\cdot,\cdot]_{\varphi},\varphi)$ be a multiplicative Hom-Lie algebra with the induced Lie algebra $(\L,[\cdot,\cdot])$.\\
(1) If $(V,\rho_{\phi},\phi)$ be a representation of  $(\L,[\cdot,\cdot]_{\varphi},\varphi)$, where $\phi$ is invertible. Then $(V,\rho)=(V,{\phi}^{-1}\circ\rho_{\phi})$ is a representation of $(\L,[\cdot,\cdot])$.\\
(2) Suppose that $(V,\rho)$ is a representation of $(\L,[\cdot,\cdot])$. If there exists $\phi\in End(V)$ such that (\ref{homrep2}), then $(V,\rho_{\phi},\phi)$ is a representation of  $(\L,[\cdot,\cdot]_{\varphi},\varphi)$.
\end{theo}

From Theorem \ref{rep-homlielie}, we have known the relationship between the representations of a multiplicative Hom-Lie algebra and the representations of its induced Lie algebra. Similar with \cite{[ABM]}, the following theorem gives the closely relation between the irreducible representations of Lie algebra and the irreducible representations of the corresponding multiplicative Hom-Lie algebra.
\begin{theo}\label{proprep-homlielie}
Let $(\L,[\cdot,\cdot]_{\varphi},\varphi)$ be a multiplicative Hom-Lie algebra and $(V,\rho_{\phi},\phi)$ a representation with $\phi$ invertible. If $(V,\rho)$ is an irreducible representation of the induced Lie algebra. Then $(V,\rho_{\phi}:=\phi\circ\rho,\phi)$ is an irreducible representation of the multiplicative  Hom-Lie algebra.
\end{theo}
\noindent{\it Proof.~}
Suppose $(V,\rho_{\phi},\phi)$ is reducible. Then, there exists $W\neq{0}$ a subspace of $V$ such that $(V,\rho_{\phi},\phi|W)$ is a submodule of $(\L,\rho_{\phi},\phi)$. That is $\phi(W)\subset{W}$ and $\rho_{\phi}(x)W\subset{W}$, $\forall x\in\L$. Hence, $\phi\circ\rho(x)W\subset{W}$,  $\forall x\in\L$ and then $\rho(x)W\subset{\phi^{-1}(W)}\subset{W}$, $\forall x\in\L$. So $W$ is submodule for $(V,\rho)$. This is a contradiction.
\hfill$\Box$

By Definition \ref{heivir} and Proposition \ref{yautwist}, we obtain the twisted Heisenberg-Virasoro algebra of Hom type, which is a multiplicative Hom-Lie algebra. In \cite{[LTL]}, the authors proved that there are nontrivial Hom-Lie algebra structure on the twisted Heisenberg-Virasoro algebra.
\begin{lemm}
\label{homstr}
As for the twisted Heisenberg-Virasoro algebra $\L$, its nontrivial Hom-Lie algebra stucture $(\L,[\cdot,\cdot]_{\varphi},\varphi)$ is given as follows:
\begin{align}
\label{b1}&\varphi(L_n)=L_n+dI_n+\delta_{n,0}dC_{LI}-\delta_{n,0}\frac{d^2}{2}C_I,\\
\label{b2}&\varphi(I_n)=I_n-\delta_{n,0}dC_I,\\
\label{b3}&\varphi(C_L)=C_L,\\
\label{b4}&\varphi(C_{LI})=C_{LI},\\
\label{b5}&\varphi(C_I)=C_I,
\end{align}
where $\varphi:\L\rightarrow\L, d\in\C$, $\forall n\in\Z$.
\end{lemm}
\vskip7pt

\cl{\bf\S3. \ Irreducible representations of the twisted Heisenberg-Virasoro algebra of Hom type}
\setcounter{section}{3}\setcounter{theo}{0}\setcounter{equation}{0}
By Theorem \ref{hcmod}, we know that the Harish-Chandra modules of the intermediate series over the twisted Heisenberg-Virasoro algebra have seven class: $A(\a,\b,F)$, $A(\a,F)$, $B(\a,F)$, $U_F$, $V_F$, $\widetilde{U}_F$, $\widetilde{V}_F$.
In this section, we suppose the seven class Harish-Chandra modules of the intermediate series over the twisted Heisenberg-Virasoro algebra are irreducible. So Thus, using Theorem \ref{rep-homlielie}, \ref{proprep-homlielie} and  Lemma \ref{homstr}, we can
use the irreducible representations of the twisted Heisenberg-Virasoro algebra to construct the irreducible representations of the twisted Heisenberg-Virasoro algebra of Hom-type.

We can suppose these modules of Hom type are $A^{\prime}(\a,\b,F)$, $A^{\prime}(\a,F)$, $B^{\prime}(\a,F)$, $U^{\prime}_F$, $V^{\prime}_F$, $\widetilde{U}^{\prime}_F$, $\widetilde{V}^{\prime}_F$. We need to note, if the linear map $\phi\in End(V)$ is invertible, these representations of Hom type are also irreducible representations, and these modules of Hom type are not always weight modules.

Next, we give irreducible representations over the twisted Heisenberg-Virasoro algebra of Hom type respectively by Theorem \ref{proprep-homlielie}.
By Lemma \ref{homstr}, we know that there are nontrivial Hom-Lie algebra structures on the twisted Heisenberg-Virasoro algebra.
Then the irreducible representations of the twisted Heisenberg-Virasoro algebra of Hom-type $\L_{\varphi}$ can be induced by the irreducible representation of its induced Lie algebra $\L$.
\begin{lemm}
\label{aabf}
Let $A(\alpha, \beta, F)$ be the irreducible Harish-Chandra modules over $\L$ defined by (\ref{a1}). Then the irreducible representations $A'(\alpha, \beta, F)$ of $\L_{\varphi}$ induced by $A(\alpha, \beta, F)$ can be defined as follows.
$$L_n(v_t)=(\a+t+{\b}n)a^{n+t}mv_{k(n+t)+q}, I_n(v_t)=Fa^{n+t}mv_{k(n+t)+q},  $$
where $m=a_{0,q}$, $q=k\a-\a-kFd\in \mathbb{Z}$.
\end{lemm}
\noindent{\it Proof.~}
(1) $A(\a,\b,F)$:\\
by $\rho((\varphi(x))\circ\phi=\phi\circ\rho(x)$,
we can get,
when $x=L_n$,
\begin{center}
$\phi\circ\rho(L_n)=\rho((\varphi(L_n))\circ\phi, n\in\Z$, \end{center}
then \begin{center} $\phi\circ\rho(L_n)(v_t)=\rho((\varphi(L_n))\circ\phi(v_t), n,t\in\Z$.
\end{center}
Suppose $\phi(v_t)=\sum_{j\in\Z}a_{t,j}v_j$, by (2.1), then
\begin{center}
$\phi((L_n)v_{n+t})=(\frac{1}{k}a^nL_{kn}+a^n(cn+d)I_{kn})(\sum_{j\in\Z}a_{t,j}v_j), n,t\in\Z$,\\
\end{center}
since $L_n(v_t)=(\a+t+{\b}n)v_{n+t}$, $I_n(v_t)=Fv_{n+t}$, then
\begin{center}
$(\a+t+{\b}n)\sum_{j\in\Z}a_{n+t,n+j}v_{n+j}=\frac{1}{k}a^n\sum_{j\in\Z}a_{t,j}(a+j+{\b}kn)v_{kn+j}+Fa^n(cn+d)\sum_{j\in\Z}a_{t,j}v_{kn+j}, n,t\in\Z$.
\end{center}
By comparing the coefficients, we have
\begin{align}
(\a+t+{\b}n)a_{n+t,n+j}=(\frac{1}{k}a^n(\a+(-((k-1)n-j))+{\b}kn)+Fa^n(cn+d))a_{t,(-((k-1)n-j)},\label{3.1}
\end{align}
where $n,t,j\in\Z$.\\
When $x=I_n$,
\begin{center}
$\phi\circ\rho(I_n)(v_t)=\rho((\varphi(I_n))\circ\phi(v_t), n,t\in\Z$,
\end{center}
similarly, we have
\begin{center}
$\phi(Fv_{n+t})=a^nbI_{kn}(\sum_{j\in\Z}a_{t,j}v_j), n,t\in\Z$,
\end{center}
then
\begin{center}
$F\sum_{j\in\Z}a_{n+t,n+j}v_{n+j}=Fa^nb\sum_{j\in\Z}a_{t,j}v_{kn+j}, n,t\in\Z$,
\end{center}
and
\begin{align}
Fa_{n+t,n+j}=Fa^nba_{t,(-((k-1)n-j)}, n,t,j\in\Z.\label{3.2}
\end{align}
Next, we consider four cases:\\
(i)$n=0, t=0$;  (ii)$n=0, t\neq0$; (iii)$n\neq0, t=0$;  (iv)$n\neq0, t\neq0$.\\
Case (i) $n=0, t=0$:\\
 by (\ref{3.1}) and (\ref{3.2}) (each of the following cases uses these two formulas), we have
\begin{align*}
\a{a_{0,j}}&=(\frac{1}{k}(\a+j)+Fd)a_{0,j}, \text{$j\in\Z$,}\\
Fa_{0,j}&=Fba_{0,j}, \text{$j\in\Z$.}
\end{align*}
Thus
\begin{align*}
a_{0,j}=0,& \text{ when $j\neq{k\a}-{\a}-kFd\in\Z$,}\\
a_{0,j}\neq0,& \text{ when $b=1,j\in\Z$.}
\end{align*}
Case (ii) $n=0, t\neq0$:
\begin{align*}
(\a+t){a_{t,j}}&=(\frac{1}{k}(\a+j)+Fd)a_{t,j}, \text{$t,j\in\Z$,}\\
Fa_{t,j}&=Fba_{t,j}, \text{$t,j\in\Z$.}
\end{align*}
Thus
\begin{align*}
a_{t,j}=0,& \text{ when $j\neq{k\a}-{\a}-kFd+kt\in\Z,t\in\Z$,}\\
a_{t,j}\neq0,& \text{ when $b=1,j\in\Z,t\in\Z$.}
\end{align*}
Case (iii) $n\neq0, t=0$:
\begin{align*}
(\a+{\b}n){a_{n,n+j}}&=(\frac{1}{k}a^n(\a+(-((k-1)n-j))+{\b}kn)+Fa^n(cn+d))a_{0,-((k-1)n-j)}, \text{$n,j\in\Z$,}\\
a_{n,n+j}&=a^na_{0,-((k-1)n-j)}, \text{$n\in\Z$.}
\end{align*} By these two equations,we have\\ $(\a+{\b}n)a^na_{0,-((k-1)n-j)}=(\frac{1}{k}a^n(\a+(-((k-1)n-j))+{\b}kn)+Fa^n(cn+d))a_{0,-((k-1)n-j)},n,j\in\Z$.
Then $$((\frac{1}{k}-1)\a+(\frac{1}{k}-1+Fc)+\frac{1}{k}j+Fcn+Fd)a_{0,-((k-1)n-j)}=0,n,j\in\Z.$$\\
since $\frac{1}{k}\neq0$, we have $a_{0,-((k-1)n-j)}=0$.
In particular, when $-((k-1)n-j)={k\a}-{\a}-kFd$, we have \\
$$(\a+{\b}n)a_{0,{k\a}-{\a}-kFd}=(\a+{\b}n+Fcn)a_{0,{k\a}-{\a}-kFd},n\in\Z,$$\\
if $a_{0,k\a-\a-kFd}\neq0$, then $c=0$.\\
Case (iv) $n\neq0, t\neq0$ (similarly to case (iii)):
\begin{align*}
(\a+t+{\b}n)a_{n+t,n+j}&=(\frac{1}{k}a^n(\a+(-((k-1)n-j))+{\b}kn)+Fa^n(cn+d))a_{t,(-((k-1)n-j)}, \text{$n,t,j\in\Z$,}\\
Fa_{n+t,n+j}&=Fa^na_{t,(-((k-1)n-j)}, \text{$n,t,j\in\Z$.}
\end{align*}
Similarly, since $\frac{1}{k}\neq0$, we have $a_{t,-((k-1)n-j)}=0$, and $a_{t,k\a-\a-kFd}\neq0$, when $c=0$.

To sum up,
when the parameters of nonzero endomorphism of the twisted Heisenberg-Virasoro algebra $\varphi$ satisfy\\
\begin{equation*}
\begin{cases}
b=1,\\
c=0,\\
k\a-\a-kFd\in\Z,
\end{cases}
\end{equation*}
we have
\begin{center}
$a_{t,k\a-\a-kFd+kn}=a^ta_{0,k\a-\a-kFd},n,t\in\Z$.
\end{center}
So $\phi(v_t)=a_{t,tk+q}v_{tk+q}=a^tmv_{tk+q}$, where $m=a_{0,q}$, $q=k\a-\a-kFd$.

Hence $A^{\prime}(\a,\b,F)$:
$L_n(v_t)=(\a+t+{\b}n)a^{n+t}mv_{k(n+t)+q}. I_n(v_t)=Fa^{n+t}mv_{k(n+t)+q}$, where $m=a_{0,q}$, $q=k\a-\a-kFd$.

\begin{lemm}
\label{abf}
Let $A(\a,F)$ be the irreducible Harish-Chandra modules over $\L$ defined by (\ref{a2}). Then the irreducible representations $A^{\prime}(\a,F)$ of $\L_{\varphi}$ induced by $A(\alpha, F)$ can be defined as follows.
$L_n(v_t)=(t+n)ka^{n+t}m_{1}v_{k(n+t)}, t\neq0, L_n(v_0)=n(n+\a)a^nm_{1}v_{kn}. I_n(v_t)=0, t\neq0, I_n(v_0)=nFa^nm_{1}v_{kn}$, where $m_{1}=a_{0,0}$.
\end{lemm}
\noindent{\it Proof.~}

(2) $A(\a,F)$:\\
when $x=L_n$,
\begin{center}
$\phi\circ\rho(L_n)(v_t)=\rho((\varphi(L_n))\circ\phi(v_t),n,t\in\Z$,\\
\end{center}
because of $L_n(v_t)=(t+n)v_{n+t}, t\neq0$. Then\\
\begin{center}
$\phi((t+n)v_{n+t})=(\frac{1}{k}a^nL_{kn}+a^n(cn+d)I_{kn})(\sum_{j\in\Z}a_{t,j}v_j),n,t\in\Z$,\\
\end{center}
since $I_n(v_t)=0, t\neq0, I_n(v_0)=nFv_n$, then
\begin{align*}
(t+n)\sum_{j\in\Z}a_{n+t,n+j}v_{n+j}&=\frac{1}{k}a^n\sum_{j\neq0,j\in\Z}a_{t,j}(kn+j)v_{kn+j}+\frac{1}{k}a^na_{t,0}kn(kn+\a)v_{kn}\\
&+knFa^n(cn+d)a_{t,0}v_{kn}, n,t\in\Z.
\end{align*}
Since the coefficients of $v_{n+j}$ and $v_{kn}$ are equal on both sides respectively, we have
\begin{align}
(t+n)a_{n+t,n+j}&=\frac{1}{k}a^n(n+j)a_{t,-((k-1)n-j)}, \text{$n,t,j\in\Z$},\label{3.3}\\
(t+n)a_{n+t,kn}&=(\frac{1}{k}a^nkn(kn+\a)+knFa^n(cn+d))a_{t,0}, \text{$n,t,\in\Z$}.\label{3.4}
\end{align}
when $t=0$,
\begin{center}
$\phi\circ\rho(L_n)(v_0)=\rho((\varphi(L_n))\circ\phi(v_0),n\in\Z$,\\
\end{center}
because of $L_n(v_0)=n(n+\a)v_n$, then
\begin{center}
$\phi(n(n+\a)v_n)=(\frac{1}{k}a^nL_{kn}+a^n(cn+d)I_{kn})(\sum_{j\in\Z}a_{0,j}v_j),n\in\Z$,\\
\end{center}
and
\begin{align*}
(n+\a)\sum_{j\in\Z}a_{n,j}v_j&=\frac{1}{k}a^n\sum_{j\neq0,j\in\Z}a_{0,j}(kn+j)v_{kn+j}+\frac{1}{k}a^na_{0,0}kn(kn+\a)v_{kn}\\
&+knFa^n(cn+d)a_{0,0}v_{kn},n\in\Z.
\end{align*}
Since the coefficients of $v_j$ and $v_{kn}$ are equal on both sides respectively, we have
\begin{align}
n(n+\a)a_{n,j}&=\frac{1}{k}a^nja_{0,-kn+j}, \text{$n,j\in\Z$,} \label{3.5}\\
n(n+\a)a_{n,kn}&=(\frac{1}{k}a^nkn(kn+\a)+knFa^n(cn+d))a_{0,0}, \text{$n\in\Z$,} \label{3.6}
\end{align}
When $x=I_n$,
\begin{center}
$\phi\circ\rho(I_n)(v_t)=\rho((\varphi(I_n))\circ\phi(v_t),n,t\in\Z$,\\
\end{center}
because of $I_n(v_t)=0, t\neq0$, then
\begin{center}
$a^nbI_{kn}(\sum_{j\in\Z}a_{t,j}v_j)=0,n,t\in\Z$.
\end{center}
then
\begin{align}
knFa^nba_{t,0}=0,n,t\in\Z.\label{3.7}
\end{align}
When $t=0$,
\begin{center}
$\phi\circ\rho(I_n)(v_0)=\rho((\varphi(I_n))\circ\phi(v_0),n\in\Z$,\\
\end{center}
because of $I_n(v_0)=nFv_n$, then
\begin{center}
$nF\sum_{j\in\Z}a_{n,j}v_j=knFa^nba_{0,0}v_{kn},n\in\Z$.
\end{center}
By comparing the coefficients of $v_j$($j\neq{kn}\in\Z$) and $v_{kn}$, we have
\begin{align}
a_{n,j}&=0, \text{ when $j\neq{kn}\in\Z, n\neq0\in\Z$,} \label{3.8}\\
a_{n,kn}&=ka^nba_{0,0}, \text{ when $n\neq0\in\Z,b\neq0$.} \label{3.9}
\end{align}
Case (i) $n=0, t=0$:\\
we substitute case(i) into (\ref{3.5}) and (\ref{3.6}),
\begin{center}
$\frac{1}{k}ja_{0,j}=0,j\in\Z$.
\end{center}
So $a_{0,j}=0$, when $j\neq0\in\Z$.\\
Case (ii) $n=0, t\neq0$:\\
by (\ref{3.3}), (\ref{3.4}) and (\ref{3.7}), we have
\begin{align*}
ta_{t,j}&=\frac{1}{k}ja_{t,j},\text{$t,j\in\Z$,}\\
ta_{t,0}&=0,\text{$t\in\Z$.}
\end{align*}
So $a_{t,j}=0$, when $j\neq{tk}\in\Z$.\\
Case (iii) $n\neq0, t=0$:\\
by (\ref{3.5}), (\ref{3.6}), (\ref{3.8}) and (\ref{3.9}), we have
\begin{align*}
n(n+\a)a_{n,kn}&=(\frac{1}{k}a^nkn(kn+\a)+knFa^n(cn+d))a_{0,0},n\in\Z,\\
a_{n,j}&=0,\text{ when $j\neq{kn}\in\Z, n\neq0\in\Z$,}\\
a_{0,j}&=0, \text{ when $j\neq0\in\Z$,}\\
a_{n,kn}&=ka^nba_{0,0},\text{ when $n\neq0\in\Z,b\neq0$.}
\end{align*}
So $\a\notin\Z/\{0\}$.\\
Case (iv) $n\neq0, t\neq0, t=-n$:\\
by (\ref{3.3}), (\ref{3.4}) and (\ref{3.7}) we have
$$a_{t,0}=0,t\in\Z.$$
Case (v) $n\neq0, t\neq0, t\neq-n$:\\
by (\ref{3.3}), (\ref{3.4}) and (\ref{3.7}) we have
\begin{align*}
(t+n)a_{n+t,n+j}&=\frac{1}{k}a^n(n+j)a_{t,-((k-1)n-j)},n,t,j\in\Z,\\
a_{n+t,kn}&=0,n,t\in\Z,\\
a_{t,0}&=0,t\in\Z.
\end{align*}
In particular,
\begin{align}
a_{n+t,k(n+t)}=a^na_{t,tk},n,t\in\Z.\label{3.10}
\end{align}
Through (\ref{3.6}), (\ref{3.9}) and (\ref{3.10}), we can get
\begin{align*}
n(n+\a)a_{n,kn}&=(\frac{1}{k}a^nkn(kn+\a)+knFa^n(cn+d))a_{0,0},n\in\Z,\\
a_{n+t,k(n+t)}&=a^na_{t,tk},n,t\in\Z,\\
a_{n,kn}&=ka^nba_{0,0},n\in\Z.
\end{align*}
Thus
$(1-k-kFc)n+(-kFd)a_{0,0}=0$, if $a_{0,0}\neq0$, we have to make $kb=1,
1-k-kFc=0,
d=0$.

To sum up,
when the parameters of nonzero endomorphism of the twisted Heisenberg-Virasoro algebra $\varphi$ satisfy\\
\begin{equation*}
\begin{cases}
\a\notin\Z/\{0\},\\
kb=1,\\
1-k-kFc=0,\\
d=0,
\end{cases}
\end{equation*}
we get
\begin{center}
$a_{t,kt}=a^ta_{0,0},t\in\Z$.
\end{center}
So $\phi(v_t)=a_{t,tk}v_{tk}=a^tm_{1}v_{tk}$, where $m_{1}=a_{0,0}$.

Hence $A^{\prime}(\a,F)$:
$L_n(v_t)=(t+n)ka^{n+t}m_{1}v_{k(n+t)}, t\neq0, L_n(v_0)=n(n+\a)a^nm_{1}v_{kn}. I_n(v_t)=0, t\neq0, I_n(v_0)=nFa^nm_{1}v_{kn}$, where $m_{1}=a_{0,0}$.\\

\begin{lemm}
\label{baf}
Let $B(\a,F)$ be the irreducible Harish-Chandra modules over $\L$ defined by (\ref{a3}). Then the irreducible representations $A^{\prime}(\a,F)$ of $\L_{\varphi}$ induced by $A(\alpha, F)$ can be defined as follows.
$L_n(v_t)=(t+n)ka^{n+t}m_{1}v_{k(n+t)}, t\neq0, L_n(v_0)=n(n+\a)a^nm_{1}v_{kn}. I_n(v_t)=0, t\neq0, I_n(v_0)=nFa^nm_{1}v_{kn}$, where $m_{1}=a_{0,0}$.
\end{lemm}
\noindent{\it Proof.~}

(3) $B(\a,F)$:\\
when $x=L_n$,
\begin{center}
$\phi\circ\rho(L_n)(v_t)=\rho((\varphi(L_n))\circ\phi(v_t), n,t\in\Z$,\\
\end{center}
because of $L_n(v_t)=tv_{n+t}, t\neq{-n}$, then
\begin{center}
$\phi(tv_{n+t})=(\frac{1}{k}a^nL_{kn}+a^n(cn+d)I_{kn})(\sum_{j\in\Z}a_{t,j}v_j),n,t\in\Z$,\\
\end{center}
since $I_n(v_t)=0, t\neq{-n}, I_n(v_{-n})=nFv_0$, then
\begin{align*}
t\sum_{j\in\Z}a_{n+t,n+j}v_{n+j}&=\frac{1}{k}a^n\sum_{j\neq{-kn},j\in\Z}a_{t,j}jv_{kn+j}+\frac{1}{k}a^na_{t,-kn}(-kn)(kn+\a)v_0\\
&+knFa^n(cn+d)a_{t,-kn}v_0,n,t\in\Z.
\end{align*}
Since the coefficients of $v_{n+j}$ and $v_0$ are equal on both sides respectively, we have
\begin{align*}
ta_{n+t,n+j}&=\frac{1}{k}a^n(-((k-1)n-j))a_{t,-((k-1)n-j)},\text{$n,t,j\in\Z$,}\\
ta_{n+t,0}&=(\frac{1}{k}a^n(-kn)(kn+\a)+knFa^n(cn+d))a_{t,-kn},\text{$n,t\in\Z$,}
\end{align*}
when $t=-n$,
\begin{center}
$\phi\circ\rho(L_n)(v_{-n})=\rho((\varphi(L_n))\circ\phi(v_{-n}),n\in\Z$,\\
\end{center}
because of $L_n(v_{-n})=-n(n+\a)v_0$, then
\begin{center}
$\phi(-n(n+\a)v_0)=(\frac{1}{k}a^nL_{kn}+a^n(cn+d)I_{kn})(\sum_{j\in\Z}a_{-n,j}v_j),n\in\Z$,\\
\end{center}
since $I_n(v_t)=0, t\neq{-n}, I_n(v_{-n})=nFv_0$, then
\begin{align*}
-n(n+\a)\sum_{j\in\Z}a_{0,j}v_{n+j}&=\frac{1}{k}a^n\sum_{j\neq{-kn},j\in\Z}a_{-n,j}jv_{kn+j}+\frac{1}{k}a^na_{-n,-kn}(-kn)(kn+\a)v_0\\
&+knFa^n(cn+d)a_{-n,-kn}v_0,n\in\Z.
\end{align*}
Since the coefficients of $v_{n+j}$ and $v_0$ are equal on both sides respectively, we have
\begin{align*}
-n(n+\a)a_{0,j}&=\frac{1}{k}a^n(-(kn-j))a_{-n,-(kn-j)},\text{$n,j\in\Z$,}\\
-n(n+\a)a_{0,0}&=(\frac{1}{k}a^n(-kn)(kn+\a)+knFa^n(cn+d))a_{-n,-kn},\text{$n\in\Z$.}
\end{align*}
When $x=I_n$,
\begin{center}
$\phi\circ\rho(I_n)(v_t)=\rho((\varphi(I_n))\circ\phi(v_t),n,t\in\Z$,\\
\end{center}
because of $I_n(v_t)=0, t\neq{-n}$, then
\begin{center}
$a^nbI_{kn}(\sum_{j\in\Z}a_{t,j}v_j)=0,n,t\in\Z$.
\end{center}
then
\begin{align*}
knFa^nba_{t,-kn}=0,n,t\in\Z,
\end{align*}
when $t=-n$,
\begin{center}
$\phi\circ\rho(I_n)(v_{-n})=\rho((\varphi(I_n))\circ\phi(v_{-n}),n\in\Z$,\\
\end{center}
because of $I_n(v_{-n})=nFv_0$, then
\begin{center}
$nF\sum_{j\in\Z}a_{0,j}v_j=knFa^nba_{-n,-kn}v_0,n\in\Z$.
\end{center}
By comparing the coefficients of $v_j$($j\neq0$) and $v_0$, we have
\begin{align*}
a_{0,j}&=0, \text{ when $j\neq0\in\Z, n\neq0\in\Z$,}\\
a_{0,0}&=ka^nba_{-n,-kn}, \text{ when $n\neq0\in\Z,b\neq0$.}
\end{align*}
Next, we still consider five cases (i) $n=0, t=0$, (ii) $n=0, t\neq0$, (iii) $n\neq0, t=0$, (iv) $n\neq0, t\neq=0, t\neq-n$, (v) $n\neq0, t\neq=0, t=-n$.\\

Similarly, to sum up,
when the parameters of nonzero endomorphism $\varphi$ satisfy\\
\begin{equation*}
\begin{cases}
\a\notin\Z/\{0\},\\
k^{-1}b^{-1}=1,\\
1-k+kFc=0,\\
d=0,
\end{cases}
\end{equation*}
we have
\begin{center}
$a_{t,kt}=a^ta_{0,0},t\in\Z$.
\end{center}
So $\phi(v_t)=a_{t,tk}v_{tk}=a^tm_{1}v_{tk}$, where $m_{1}=a_{0,0}$.

Hence $B^{\prime}(\a,F)$:
$L_n(v_t)=ta^{n+t}m_{1}v_{k(n+t)}, t\neq-n, L_n(v_{-n})=-n(n+\a)a^nm_{1}v_{0}. I_n(v_t)=0, t\neq-n, I_n(v_{-n})=nFm_{1}v_0$, where $m_{1}=a_{0,0}$.\\

\begin{lemm}
\label{baf}
Let $U_F$ be the irreducible Harish-Chandra modules over $\L$ defined by (\ref{a4}). Then the irreducible representations $U_F'$ of $\L_{\varphi}$ induced by $U_F$ can be defined as follows.
$L_n(v_t)=ta^{n+t}m_{1}v_{k(n+t)}, t\neq-n. I_n(v_t)=0, t\neq-n, I_n(v_{-n})=nFm_{1}v_0$, where $m_{1}=a_{0,0}$.
\end{lemm}
\noindent{\it Proof.~}
(4) $U_F$(similar to $B(\a,F)$):\\
when $x=L_n$,
Similarly, since the coefficients of $v_{n+j}$ and $v_0$ are equal on both sides respectively, we have
\begin{align*}
ta_{n+t,n+j}&=\frac{1}{k}a^n(-((k-1)n-j))a_{t,-((k-1)n-j)},\text{$n,t,j\in\Z$,}\\
ta_{n+t,0}&=(\frac{1}{k}a^n(-kn)+knFa^n(cn+d))a_{t,-kn},\text{$n,t\in\Z$.}
\end{align*}
When $x=I_n$, similarly,\\
$t\neq-n$:
\begin{align*}
knFa^nba_{t,-kn}=0,n,t\in\Z.
\end{align*}
$t=-n$:
\begin{align*}
a_{0,j}&=0, \text{ when $j\neq0\in\Z, n\neq0\in\Z$,}\\
a_{0,0}&=ka^nba_{-n,-kn}, \text{ when $n\neq0\in\Z,b\neq0$.}
\end{align*}
Similarly, we consider the values of $n$ and $t$, when\\
\begin{equation*}
\begin{cases}
k^{-1}b^{-1}=1,\\
c=0,\\
d=0,
\end{cases}
\end{equation*}
we have
\begin{center}
$a_{t,kt}=a^ta_{0,0},t\in\Z$.
\end{center}
So $\phi(v_t)=a_{t,tk}v_{tk}=a^tm_{1}v_{tk}$, where $m_{1}=a_{0,0}$.

Hence $U^{\prime}_F$:
$L_n(v_t)=ta^{n+t}m_{1}v_{k(n+t)}, t\neq-n. I_n(v_t)=0, t\neq-n, I_n(v_{-n})=nFm_{1}v_0$, where $m_{1}=a_{0,0}$.\\

\begin{lemm}
\label{vf}
Let $V_F$ be the irreducible Harish-Chandra modules over $\L$ defined by (\ref{a5}). Then the irreducible representations $V_F'$ of $\L_{\varphi}$ induced by $U_F$ can be defined as follows:
$L_n(v_t)=(t+n)a^{n+t}m_{1}v_{k(n+t)}, t\neq-n. I_n(v_t)=0, t\neq0, I_n(v_0)=nFa^nm_{1}v_kn$, where $m_{1}=a_{0,0}$.
\end{lemm}

(5) $V_F$(similar to $A(\a,F)$):\\
when $x=L_n$,
Similarly, since the coefficients of $v_{n+j}$ and $v_{kn}$ are equal on both sides respectively, we have
\begin{align*}
(t+n)a_{n+t,n+j}&=\frac{1}{k}a^n(n+j))a_{t,-((k-1)n-j)},\text{$n,t,j\in\Z$,}\\
(t+n)a_{n+t,kn}&=(\frac{1}{k}a^n(kn)+knFa^n(cn+d))a_{t,0},\text{$n,t\in\Z$.}
\end{align*}
When $x=I_n$, similarly,\\
$t\neq0$:
\begin{align*}
knFa^nba_{t,0}=0,n,t\in\Z.
\end{align*}
$t=0$:
\begin{align*}
a_{n,j}&=0, \text{ when $j\neq{kn}\in\Z, n\neq0\in\Z$,}\\
a_{n,kn}&=ka^nba_{0,0}, \text{ when $n\neq0\in\Z,b\neq0$.}
\end{align*}
Similarly, when
\begin{equation*}
\begin{cases}
kb=1,\\
c=0,\\
d=0,
\end{cases}
\end{equation*}
we have
\begin{center}
$a_{t,kt}=a^ta_{0,0},t\in\Z$.
\end{center}
So $\phi(v_t)=a_{t,tk}v_{tk}=a^tm_{1}v_{tk}$, where $m_{1}=a_{0,0}$.

Hence $V^{\prime}_F$:
$L_n(v_t)=(t+n)a^{n+t}m_{1}v_{k(n+t)}, t\neq-n. I_n(v_t)=0, t\neq0, I_n(v_0)=nFa^nm_{1}v_kn$, where $m_{1}=a_{0,0}$.\\

\begin{lemm}
\label{uf}
Let $\widetilde{U}_F$ be the irreducible Harish-Chandra modules over $\L$ defined by (\ref{a6}). Then the irreducible representations $\widetilde{U}_F$ of $\L_{\varphi}$ induced by $U_F$ can be defined as follows:
$L_n(v_t)=(t+n)a^{n+t}m_{1}v_{k(n+t)}, t\neq-n, L_n(v_{-n})=0. I_n(v_t)=0, t\neq-n, I_n(v_{-n})=nFm_{1}v_0$, where $m_{1}=a_{0,0}$.
\end{lemm}

(6) $\widetilde{U}_F$:\\
when $x=L_n$,
since $\phi\circ\rho(L_n)(v_t)=\rho((\varphi(L_n))\circ\phi(v_t),n,t\in\Z$,\\
$t\neq-n$:\\
Similarly, we can get $$t\sum_{j\in\Z}a_{n+t,n+j}v_{n+j}=\frac{1}{k}a^n\sum_{j\neq{-kn},j\in\Z}a_{t,j}jv_{kn+j}+knFa^n(cn+d)a_{t,-kn}v_0,n,t\in\Z$$ since the coefficients of $v_{n+j}$ and $v_0$ are equal on both sides respectively, we have
\begin{align}
ta_{n+t,n+j}&=\frac{1}{k}a^n(-((k-1)n-j))a_{t,-((k-1)n-j)},\text{$n,t,j\in\Z$,} \label{3.11}\\
ta_{n+t,0}&=knFa^n(cn+d))a_{t,-kn},\text{$n,t\in\Z$.} \label{3.12}
\end{align}
$t=-n$:\\
because of $L_n(v_{-n})=0$, then
\begin{align}
\cdots+\frac{1}{k}a^na_{-n,-1}(-1)v_{kn-1}+\cdots+\frac{1}{k}a^na_{-n,j}(j)v_{kn+j}+\cdots
+knFa^n(cn+d)a_{-n,-kn}v_0=0,n,j\in\Z, \label{3.13}
\end{align}
where $j\in\Z$.\\
When $x=I_n$(same as before),
since $\phi\circ\rho(I_n)(v_t)=\rho((\varphi(I_n))\circ\phi(v_t),n,t\in\Z$,\\
$t\neq-n$:
\begin{align}
knFa^nba_{t,-kn}=0,n,t\in\Z. \label{3.14}
\end{align}
$t=-n$:
\begin{align}
a_{0,j}&=0, \text{  $j\neq0\in\Z, n\neq0\in\Z$,} \label{3.15}\\
a_{0,0}&=ka^nba_{-n,-kn}, \text{ when $n\neq0\in\Z,b\neq0$.} \label{3.16}
\end{align}
Case (i) $n=0, t=0$:\\
by (\ref{3.13}), we have
\begin{center}
$\cdots, a_{0,-1}, a_{0,1}, \cdots, a_{0,j}=0, \cdots,j\in\Z$.
\end{center}
Case (ii) $n=0, t\neq0$:\\
by (\ref{3.11}), (\ref{3.12}) and (\ref{3.14}), we have
\begin{align*}
ta_{t,j}&=\frac{1}{k}ja_{t,j},\text{$t,j\in\Z$,}\\
ta_{t,0}&=0,\text{$t\in\Z$.}
\end{align*}
So $a_{t,j}=0$, when $j\neq{tk}\in\Z$.\\
Case (iii) $n\neq0, t=0$:\\
by (\ref{3.11}), (\ref{3.12}) and (\ref{3.14}), we have
\begin{align*}
&\frac{1}{k}a^n(-((k-1)n-j))a_{0,-((k-1)n-j)}=0,\text{$n,j\in\Z$,}\\
&a_{0,-kn}=0,\text{$n\in\Z$.}
\end{align*}
Case (iv) $n\neq0, t\neq0,t=-n$:\\
by (\ref{3.13}), we have
$\cdots, a_{-n,-1}, a_{-n,1}, \cdots, a_{-n,j}=0, \cdots$, and $knFa^n(cn+d)a_{-n,-kn}=0$, if $a_{-n,-kn}\neq0$, then $c=0, d=0$.\\
Case (v) $n\neq0, t\neq0,t\neq-n$:\\
by (\ref{3.11}), (\ref{3.12}) and (\ref{3.14}), we have
\begin{align*}
ta_{n+t,n+j}&=\frac{1}{k}a^n(-((k-1)n-j))a_{t,-((k-1)n-j)}, n,t,j\in\Z,\\
ta_{n+t,0}&=0,n,t\in\Z,\\
a_{t,-kn}&=0,n,t\in\Z.
\end{align*}
In particular,\\
\begin{align}
a_{n+t,k(n+t)}=a^na_{t,tk},n,t\in\Z.\label{3.17}
\end{align}
Through (\ref{3.16}) and (\ref{3.17}), we can get
\begin{align*}
a_{-n,-kn}&=k^{-1}a^{-n}b^{-1}a_{0,0},n\in\Z,\\
a_{t,tk}&=a^{-n}a_{n+t,k(n+t)},n,t\in\Z.
\end{align*}
Thus $k^{-1}b^{-1}=1$.
To sum up,
when\\
\begin{equation*}
\begin{cases}
k^{-1}b^{-1}=1,\\
c=0,\\
d=0,
\end{cases}
\end{equation*}
we have
\begin{center}
$a_{t,kt}=a^ta_{0,0},t\in\Z$.
\end{center}
So $\phi(v_t)=a_{t,tk}v_{tk}=a^tm_{1}v_{tk}$, where $m_{1}=a_{0,0}$.

Hence $\widetilde{U}^{\prime}_F$:
$L_n(v_t)=(t+n)a^{n+t}m_{1}v_{k(n+t)}, t\neq-n, L_n(v_{-n})=0. I_n(v_t)=0, t\neq-n, I_n(v_{-n})=nFm_{1}v_0$, where $m_{1}=a_{0,0}$.\\

\begin{lemm}
\label{uf}
Let $\widetilde{V}_F$ be the irreducible Harish-Chandra modules over $\L$ defined by (\ref{a7}). Then the irreducible representations $\widetilde{V}_F$ of $\L_{\varphi}$ induced by $V_F$ can be defined as follows:
$L_n(v_t)=ta^{n+t}m_{1}v_{k(n+t)}, t\neq-n, L_n(v_{-n})=0. I_n(v_t)=0, t\neq0, I_n(v_0)=Fa^nm_{1}v_{kn}$, where $m_{1}=a_{0,0}$.
\end{lemm}

(7) $\widetilde{V}_F$(similar to $\widetilde{U}_F$):\\
when $x=L_n$,
since $\phi\circ\rho(L_n)(v_t)=\rho((\varphi(L_n))\circ\phi(v_t)$,\\
$t\neq-n$:\\
Similarly, we can get $t\sum_{j\in\Z}a_{n+t,n+j}v_{n+j}=\frac{1}{k}a^n\sum_{j\neq{-kn},j\in\Z}a_{t,j}jv_{kn+j}+Fa^n(cn+d)a_{t,0}v_{kn},n,t\in\Z$, since the coefficients of $v_{n+j}$, $v_0$ and $v_{kn}$ are equal on both sides respectively, we have
\begin{align}
ta_{n+t,n+j}&=\frac{1}{k}a^n(-((k-1)n-j))a_{t,-((k-1)n-j)},\text{$n,t,j\in\Z$,} \label{3.18}\\
ta_{n+t,0}&=0,n,t\in\Z, \label{3.19}\\
ta_{n+t,kn}&=Fa^n(cn+d))a_{t,0},n,t\in\Z.\label{3.20}
\end{align}
$t=-n$:\\
because of $L_n(v_{-n})=0$, then
$$\cdots+\frac{1}{k}a^na_{-n,-1}(-1)v_{kn-1}+(\frac{1}{k}a^n(0)+Fa^n(cn+d))a_{-n,0}v_{kn}+\cdots+\frac{1}{k}a^na_{-n,j}(j)v_{kn+j}+\cdots=0,n,j\in\Z$$
then
\begin{align}
\cdots, a_{-n,-1}, a_{-n,1}, \cdots, a_{-n,j}=0, \cdots,n,j\in\Z.\label{3.21}
\end{align}
When $x=I_n$, similarly,\\
$t\neq0$:
\begin{align}
Fa^nba_{t,0}=0,n,t\in\Z.\label{3.22}
\end{align}
$t=0$:
\begin{align}
a_{n,j}&=0, \text{ when $j\neq{kn}\in\Z, n\neq0\in\Z$,} \label{3.23}\\
a_{n,kn}&=a^nba_{0,0}, \text{ when $n\neq0\in\Z,b\neq0$.}\label{3.24}
\end{align}
Case (i) $n=0, t=0$:\\
by (\ref{3.21}), we have
$\cdots, a_{0,-1}, a_{0,1}, \cdots, a_{0,j}=0, \cdots$, and $Fda_{0,0}=0, a_{0,0}=ba_{0,0}$, if $a_{0,0}\neq0$, then $d=0, b=1$.\\
Case (ii) $n=0, t\neq0$:\\
by (\ref{3.18}), (\ref{3.19}), (\ref{3.20}) and (\ref{3.22}),  we have
\begin{align*}
ta_{t,j}&=\frac{1}{k}ja_{t,j},\text{$t,j\in\Z$,}\\
ta_{t,0}&=0,\text{$t\in\Z$.}
\end{align*}
So $a_{t,j}=0$, when $j\neq{tk}\in\Z$.\\
Case (iii) $n\neq0, t=0$:\\
by (\ref{3.18}), (\ref{3.19}), (\ref{3.20}), (\ref{3.23}) and (\ref{3.24}), we have
\begin{align*}
&\a_{0,-((k-1)n-j)}=0,\text{$n,j\in\Z$,}\\
&Fa^n(cn+d))a_{0,0}=0,\text{$n\in\Z$,}\\
&a_{n,j}=0, \text{ when $j\neq{kn}\in\Z, n\neq0\in\Z$,}\\
&a_{n,kn}=a^nba_{0,0}, \text{ when $n\neq0\in\Z,b\neq0$.}
\end{align*}
Then $c=0$.\\
Case (iv) $n\neq0, t\neq0,t=-n$:\\
by (\ref{3.21}) and (\ref{3.22}), we have
$$\cdots, a_{-n,-1}, a_{-n,0}, a_{-n,1}, \cdots, a_{-n,j}=0, \cdots,n,j\in\Z.$$
Case (v) $n\neq0, t\neq0,t\neq-n$:\\
by (\ref{3.18}), (\ref{3.19}), (\ref{3.20}) and (\ref{3.22}), we have
\begin{align*}
ta_{n+t,n+j}&=\frac{1}{k}a^n(-((k-1)n-j))a_{t,-((k-1)n-j)}, n,t,j\in\Z,\\
a_{t,0}&=0,t\in\Z,\\
a_{n+t,kn}&=0,n,t\in\Z,\\
a_{n+t,0}&=0,n,t\in\Z.
\end{align*}
In particular,
\begin{align}
a_{n+t,k(n+t)}=a^na_{t,tk},n,t\in\Z.\label{3.25}
\end{align}
Through (\ref{3.24}) and (\ref{3.25}), we can get
\begin{align*}
a_{n,kn}&=a^na_{0,0},n\in\Z,\\
a_{n+t,k(n+t)}&=a^na_{n,kn},n,t\in\Z.
\end{align*}
To sum up,
when\\
\begin{equation*}
\begin{cases}
b=1,\\
c=0,\\
d=0,
\end{cases}
\end{equation*}
we have
\begin{center}
$a_{t,kt}=a^ta_{0,0},t\in\Z$.
\end{center}
So $\phi(v_t)=a_{t,tk}v_{tk}=a^tm_{1}v_{tk}$, where $m_{1}=a_{0,0}$.

Hence $\widetilde{V}^{\prime}_F$:
$L_n(v_t)=ta^{n+t}m_{1}v_{k(n+t)}, t\neq-n, L_n(v_{-n})=0. I_n(v_t)=0, t\neq0, I_n(v_0)=Fa^nm_{1}v_{kn}$, where $m_{1}=a_{0,0}$.
\begin{theo}
The irreducible representations over the twisted Heisenberg-Virasoro algebra of Hom type with
$C_L=C_{LI}=C_I=0$ on a vector space $V$ equipped with a basis
$$\{\cdots, v_{-n}, \cdots, v_{0}, v_{1}, \cdots, v_{n}, \cdots\}$$
are given as follows:\\
$A^{\prime}(\a,\b,F): L_n(v_t)=(\a+t+{\b}n)a^{n+t}mv_{k(n+t)+q}.
I_n(v_t)=Fa^{n+t}mv_{k(n+t)+q}$, where $m=a_{0,q}, q=k\a-\a-kFd$.\\
$A^{\prime}(\a,F): L_n(v_t)=(t+n)ka^{n+t}m_{1}v_{k(n+t)}, t\neq0, L_n(v_0)=n(n+\a)a^nm_{1}v_{kn}.
I_n(v_t)=0, t\neq0, I_n(v_0)=nFa^nm_{1}v_{kn}$, where $m_{1}=a_{0,0}$.\\
$B^{\prime}(\a,F): L_n(v_t)=ta^{n+t}m_{1}v_{k(n+t)}, t\neq-n, L_n(v_{-n})=-n(n+\a)a^nm_{1}v_{0}.
I_n(v_t)=0, t\neq-n, I_n(v_{-n})=nFm_{1}v_0$, where $m_{1}=a_{0,0}$.\\
$U^{\prime}_F: L_n(v_t)=ta^{n+t}m_{1}v_{k(n+t)}, t\neq-n. I_n(v_t)=0, t\neq-n, I_n(v_{-n})=nFm_{1}v_0$, where $m_{1}=a_{0,0}$.\\
$V^{\prime}_F: L_n(v_t)=(t+n)a^{n+t}m_{1}v_{k(n+t)}, t\neq-n. I_n(v_t)=0, t\neq0, I_n(v_0)=nFa^nm_{1}v_kn$, where $m_{1}=a_{0,0}$.\\
$\widetilde{U}^{\prime}_F: L_n(v_t)=(t+n)a^{n+t}m_{1}v_{k(n+t)}, t\neq-n, L_n(v_{-n})=0. I_n(v_t)=0, t\neq-n, I_n(v_{-n})=nFm_{1}v_0$, where $m_{1}=a_{0,0}$.\\
$\widetilde{V}^{\prime}_F: L_n(v_t)=ta^{n+t}m_{1}v_{k(n+t)}, t\neq-n, L_n(v_{-n})=0. I_n(v_t)=0, t\neq0, I_n(v_0)=Fa^nm_{1}v_{kn}$, where $m_{1}=a_{0,0}$.\\
In particular, as for $A^{\prime}(\a,\b,F)$, when $t=kt+q$, i.e $k=1, q=0$, the vector space $V$ is a weight module. As for $A^{\prime}(\a,F)$, $B^{\prime}(\a,F)$, $U^{\prime}_F$, $V^{\prime}_F$, $\widetilde{U}^{\prime}_F$, $\widetilde{V}^{\prime}_F$, when $t=kt$, i.e $k=1$, $V$ is a weight module.
\end{theo}

\end{document}